\documentclass[10pt]{article}
\usepackage{latexsym,amssymb,amsmath,amsthm,amsfonts}
\usepackage{graphicx}
\usepackage{enumerate}
\usepackage{color}

\theoremstyle{plain}
\newtheorem{thm}{Theorem}[section]
\newtheorem{prop}[thm]{Proposition}
\newtheorem{lem}[thm]{Lemma}

\newtheorem{cor}[thm]{Corollary}
\theoremstyle{definition}
\newtheorem{definition}[thm]{Definition}

\newtheorem{example}[thm]{Example}

\numberwithin{equation}{section}

\newcommand{\f}{\frac}

\newcommand{\xx}{|x|^2}
\newcommand{\yy}{|y|^2}
\newcommand{\xy}{\langle x,y\rangle}

\newcommand{\pppp}[4]%
  {\frac{\partial^3{#1}}{\partial{#2}\partial{#3}\partial{#4}}}

\newcommand{\p}{\phi}

\newcommand{\gab}{\alpha\phi \Big(b^2,\frac{\beta}{\alpha}\Big )}
\newcommand{\pt}{\phi_2}
\newcommand{\po}{\phi_1}
\newcommand{\ptt}{\phi_{22}}
\newcommand{\pot}{\phi_{12}}

\renewcommand{\a}{\alpha}
\renewcommand{\b}{\beta}
\newcommand{\ab}{(\alpha,\beta)}

\begin{document}
\title{On general $\ab$-metrics with vanishing Douglas curvature}
\author{Hongmei Zhu\footnote{supported in  a doctoral scientific research foundation of Henan Normal University (No.5101019170130)}}
\maketitle


\begin{abstract}
In this paper, we study a class of Finsler metrics called general $\ab$-metrics, which are defined by a Riemannian metric $\a$ and a $1$-form $\b$.
We find an equation which is necessary and sufficient condition  for  such Finsler metric to be a Douglas metric. By solving this equation, we obtain all of general $\ab$-metrics with vanishing Douglas curvature under certain condition. Many new non-trivial examples are explicitly constructed.
\end{abstract}

\section{Introduction}
In Finsler geometry, one of important projective invariants is Douglas curvature, which was introduced by J. Douglas \cite{Douglas}. If two Finsler metrics $F$ and $\tilde{F}$ are projectively equivalent, then they have the same Douglas curvature. The Douglas curvature always vanishes for Riemannian metrics. Finsler metrics with vanishing Douglas curvature are called {\it Douglas metrics}. Douglas metrics form a rich class of Finsler metrics including locally projectively flat Finsler metrics.

Randers metrics are an important class of Finsler metrics, which are introduced by a physicist G. Randers in 1941. A Randers metric is of the form $F=\a+\b$, where $\a$ is a Riemannian metric and $\b$ is a 1-form. However, it can also be expressed in the following navigation form
\begin{equation*}
F = \frac{\sqrt{(1-b^2)\alpha^2+\beta^2 }}{1-b^2} +\frac{\beta}{1-b^2}.\label{Randersmetric}
\end{equation*}
It is well-known that a Randers metric is a Douglas metric if and only if $\b$ is closed for both of the above expressions \cite{Mat}
. As a generalization of Randers metrics, $\ab$-metrics are also defined by a a Riemannian metric and a 1-form and given in the form
$$
F=\a\phi(\frac{\b}{\a}),
$$
where $\phi$ is a smooth function and satisfies two additional conditions. In 2009, B. Li, Y. Shen and Z. Shen gave a characterization of Douglas $\ab$-metrics with dimension $n\geq 3$ \cite{Li-Shen-Shen}. Recently, C. Yu gave a more clear characterization. {\em If $F=\a\phi(\frac{\b}{\a})$ is a non-trivial Douglas metric, then after some special deformations, $\a$ will turn to be another Riemannian metric $\bar{\a}$ and $\b$ to be another 1-form $\bar{\b}$ such that $\bar{\b}$ is close and conformal with respect to $\bar{\a}$, i.e., $\bar{b}_{i|j}=c(x)\bar{\a}_{ij}$, where  $c(x)\neq 0$ is a scalar function on the manifold. In this case, $F$ can be reexpressed as the form $F=\bar{\a}\phi(\bar{b}^{2},\frac{\bar{\b}}{\bar{\a}})$} \cite{yct}.

In fact, many famous Douglas metrics can be also expressed  in the following form
\begin{equation}\label{generalab}
F=\gab,
\end{equation}
where $\a$ is a Riemannian metric, $\b$ is a $1$-form, $b: =\|\beta_x\|_{\alpha}$ and $\p(b^2,s)$ is a smooth function. Finsler metrics in this form are called general $\ab$-metrics \cite{yct-zhm-onan}. If $\phi=\phi(s)$ is independent of $b^2$, then $F=\a\phi(\frac{\b}{\a})$ is a $\ab$-metric. If $\a=|y|$, $\b=\xy$, then $F=|y|\phi(|x|^2,\frac{\xy}{|y|})$ is the so-called spherically symmetric Finsler metrics \cite{Mo-Tenenblat}. Moreover, general $\ab$-metrics include part of Bryant's metrics \cite{Br2, yct-zhm-onan} and part of fourth root metrics \cite{Li-Shen}.
Besides Randers metrics,
square metrics can  be expressed in the following form
\begin{equation*}\label{square}
F=\frac{(\sqrt{(1-b^2)\alpha^2+\beta^2}+\beta)^2}{(1-b^2)^2\sqrt{(1-b^2)\alpha^2+\beta^2}},
\end{equation*}
It has been shown that $F$ is a non-trivial Douglas square metrics
 if and only if
\begin{equation*}
b_{i|j} = c a_{ij},
\end{equation*}
where $c=c(x)\not=0$ is a scalar function on $M$ \cite{yct}.

In this paper, we mainly study general $\ab$-metrics with vanishing Douglas curvature. Firstly, a characterization equation for such metrics to be Douglas  metrics under a suitable condition is given (Theorem \ref{main1}). By solving this equation, we obtain all general $\ab$-metrics with vanishing Douglas curvature under certain condition(Theorem \ref{main2}).
At last, we explicitly construct some new examples (see Section 6).

Here, we will assume that $\b$ is closed and conformal with respect to $\a$, i.e. (\ref{bcij}) holds. According to the relate discussions for Douglas $\ab$-metrics  \cite{cxy-tyf, Li-Shen-Shen, matsumoto2, Mo-Tenenblat, yct},  we believe that the assumption here is reasonable and appropriate.

The main results are given below.
\begin{thm}\label{main1}
Let $F=\gab$ be a non-Riemannian general $\ab$-metric on an $n$-dimensional manifold $M$. Suppose that $\b$ satisfies
\begin{equation}
b_{i|j} = c a_{ij}, \label{bcij}
\end{equation}
where $c=c(x)\not=0$ is a scalar function on $M$ and $b_{i|j}$ is the covariant derivation of $\b$ with respect to $\a$.
Then $F$ is a Douglas metric if and only if the following  PDE holds
\begin{equation}
\ptt-2(\po-s\pot) =\{f(b^2)+g(b^2)s^2\}\{ \p-s\pt+(b^2-s^2)\ptt   \}  \label{pde1****}
\end{equation}
where $f(x)$ and $g(x)$ are two arbitrary differentiable functions.
\end{thm}
Note that $\phi_1$ means the derivation of $\phi$ with respect to the first variable $b^2$.

It should be pointed out that if the scalar function $c(x)=0$, then according to Proposition \ref{Douglas curvature}, $D^{i}{}_{jkl}=0$, namely, $F=\gab$ is a Douglas metric for any function $\phi(b^{2},s)$. So it will be regarded as a trivial case.

By solving equation (\ref{pde1****}), we have the following result

\begin{thm}\label{main2} Let $F=\gab$ be a non-Riemannian general $\ab$-metric on an $n$-dimensional manifold $M$. Suppose that $\b$ satisfies
(\ref{bcij}). Then the general solution of (\ref{pde1****}) is given by
\begin{eqnarray}\label{mainsoltion}
\phi=s\left(h(b^{2})-\int \frac{\Phi\big(\eta(b^{2},s)\big)}{s^{2}\sqrt{b^{2}-s^{2}}}ds\right),
\end{eqnarray}
where
\begin{eqnarray}\label{eta}
\eta(b^{2},s):=\frac{b^{2}-s^{2}}{e^{\int(f+gb^{2})db^{2}}-(b^{2}-s^{2})\int g e^{\int(f+gb^{2})db^{2}}db^{2}},
\end{eqnarray}
where $f$, $g$ and $h$ are arbitrary smooth functions of $b^{2}$. $\Phi$ is an  arbitrary smooth function of $\eta$. Moreover, the corresponding general $\ab$-metrics of (\ref{mainsoltion}) are of Douglas type.
\end{thm}
\textbf{Remark:} If the general $\ab$-metrics $F=\gab$ given by (\ref{mainsoltion}) are regular Finsler metrics,  then (\ref{mainsoltion}) should satisfy Lemma \ref{finsler}.

\section{Preliminaries}

Let $F$ be a Finsler metric on an $n$-dimensional manifold $M$ and $G^{i}$ be the geodesic coefficients of $F$, which are defined by
\begin{eqnarray*}
G^{i}=\frac{1}{4}g^{il}\left\{[F^{2}]_{x^{k}y^{l}}y^{k}-[F^{2}]_{x^{l}}\right\},
\end{eqnarray*}
where $(g^{ij}):=\left(\frac{1}{2}[F^{2}]_{y^{i}y^{j}}\right)^{-1}$. For a Riemannian metric, the spray coefficients are determined by its Christoffel symbols as $G^{i}(x,y)=\frac{1}{2}\Gamma^{i}_{jk}(x)y^{j}y^{k}$.

By definition, a general $(\a,\b)$-metric is given by (\ref{generalab}), where $\phi(b^{2},s)$ is a positive smooth function defined on the domain $|s|\leq b<b_o$ for some positive number (maybe infinity) $b_o$.
Then the function $F=\a\phi(b^{2},s)$ is a Finsler metric for any Riemannian metric $\alpha=\sqrt{a_{ij}(x)y^{i}y^{j}}$ and  any $1$-form $\beta=b_{i}(x)y^{i}$ if and only if $\phi(b^{2},s)$ satisfies
\begin{eqnarray}\label{ppp}
\p-s\pt>0,\quad\p-s\pt+(b^2-s^2)\ptt>0,
\end{eqnarray}
when $n\geq 3$ or
\begin{eqnarray}\label{ppp1}
\p-s\pt+(b^2-s^2)\ptt>0,
\end{eqnarray}
when $n=2$ \cite{yct-zhm-onan}.

Let $\alpha=\sqrt{a_{ij}(x)y^iy^j}$  and $\beta= b_i(x)y^i$.
Denote the coefficients of the covariant derivative of
$\b$ with respect to $\a$ by $b_{i|j}$, and let
\begin{eqnarray*}
&r_{ij}=\frac{1}{2}(b_{i|j}+b_{j|i}),~s_{ij}=\frac{1}{2}(b_{i|j}-b_{j|i}),
~r_{00}=r_{ij}y^iy^j,~s^i{}_0=a^{ij}s_{jk}y^k,&\\
&r_i=b^jr_{ji},~s_i=b^js_{ji},~r_0=r_iy^i,~s_0=s_iy^i,~r^i=a^{ij}r_j,~s^i=a^{ij}s_j,~r=b^ir_i,&
\end{eqnarray*}
where $(a^{ij}):=(a_{ij})^{-1}$ and $b^{i}:=a^{ij}b_{j}$. It is easy to see that $\b$ is closed if and only if $s_{ij}=0$.

According to \cite{yct-zhm-onan}, the spray coefficients $G^i$ of a general $(\alpha,\beta)$-metric $F=\gab$ are related to the spray coefficients ${}^\a G^i$ of
$\a$ and given by
\begin{eqnarray}\label{Gi}
G^i&=&{}^\a G^i+\a Q s^i{}_0+\left\{\Theta(-2\a Q s_0+r_{00}+2\a^2
R r)+\a\Omega(r_0+s_0)\right\}\frac{y^i}{\a}\nonumber\\
&&+\left\{\Psi(-2\a Q s_0+r_{00}+2\a^2 R
r)+\a\Pi(r_0+s_0)\right\}b^i -\a^2 R(r^i+s^i),
\end{eqnarray}
where
$$Q=\frac{\pt}{\p-s\pt},\quad R=\frac{\po}{\p-s\pt},$$
$$\Theta=\frac{(\p-s\pt)\pt-s\p\ptt}{2\p\big(\p-s\pt+(b^2-s^2)\ptt\big)},
\quad\Psi=\frac{\ptt}{2\big(\p-s\pt+(b^2-s^2)\ptt\big)},$$
$$\Pi=\frac{(\p-s\pt)\pot-s\po\ptt}{(\p-s\pt)\big(\p-s\pt+(b^2-s^2)\ptt\big)},\quad
\Omega=\frac{2\po}{\p}-\frac{s\p+(b^2-s^2)\pt}{\p}\Pi.$$

In the following, we will introduce an important projective invariant.
\begin{definition} \cite{szm-diff}
 Let
\begin{eqnarray}\label{DT}
  D^{i}_{jkl} = \frac{\partial^{3}}{\partial y^{j}\partial y^{k} \partial y^{l}}\left(G^{i}-\frac{1}{n+1}\frac{\partial G^{m}}{\partial
  y^{m}}y^{i}\right),
\end{eqnarray}
where $G^{i}$ are the spray coefficients of $F$. The tensor
$D:=D^{i}_{jkl} \partial_{i}\otimes dx^{j}\otimes dx^{k}\otimes
dx^{l}$ is called {\em Douglas tensor}. A Finsler metric is called
{\it a Douglas metric} if the Douglas tensor vanishes.
\end{definition}

We require the following result in Section 6, the proof which is omitted.
\begin{lem}\label{formu}
 Let
\begin{eqnarray}\label{In}
I_{n}:=\int s^{-2}(b^{2}-s^{2})^{\frac{n-1}{2}}ds,
\end{eqnarray}
then for any natural number $n\geq 1$, we have
\begin{enumerate}[(a)]
\item $n=2m$
\begin{eqnarray}\label{even}
 I_{2m}&=&\frac{(2m-1)!!}{(2m-2)!!}\frac{1}{s}\sum_{i=1}^{m-1}\frac{(2m-2-2i)!!}{(2m-2i+1)!!}
 \big(b^{2}\big)^{i-1}(b^{2}-s^{2})^{\frac{2m-2i+1}{2}}\nonumber\\
 &&- \frac{(2m-1)!!}{(2m-2)!!}\frac{1}{s}\big(b^{2}\big)^{m-1}\big[(b^{2}-s^{2})^{\frac{1}{2}}+s \arctan\frac{s}{\sqrt{b^{2}-s^{2}}}\big]+C_{1}.
 \end{eqnarray}
\item $n=2m+1$
\begin{eqnarray}\label{odd}
I_{2m+1}&=& \frac{(2m)!!}{(2m-1)!!}\frac{1}{s}\left[\sum_{i=1}^{m}\frac{(2m-2i-1)!!}{(2m-2i+2)!!}\big(b^{2}\big)^{i-1}(b^{2}-s^{2})^{m-i+1}-\big(b^{2}\big)^{m}\right]
+C_{2},
\end{eqnarray}
\end{enumerate}
where $C_{1}$ and $C_{2}$ are arbitrary constants.
\end{lem}

\section{Douglas curvature of general $(\alpha,\beta)$-metrics}
In this section, we will compute the Douglas curvature of a general $\ab$-metric.
\begin{prop}\label{Douglas curvature}
Let $F=\gab$ be a general $\ab$-metric on an $n$-dimensional manifold $M$. Suppose that $\b$ satisfies (\ref{bcij}), then the Douglas curvature of $F$ is given by
\begin{eqnarray}
D^{i}{}_{jkl}&=&\frac{c}{\alpha}\left\{\left[(T-sT_{2})a_{kl}+T_{22}b_{l}b_{k}\right]\delta^{i}{}_{j}+\frac{1}{\alpha^{2}}\left[\frac{s}{\alpha}(3T_{22}+sT_{222})y_{l}y_{j}-
(T_{22}+sT_{222})b_{l}y_{j}\right]b_{k}y^{i}\right\}(k\rightarrow l\rightarrow j\rightarrow k)\nonumber\\
&&-\frac{c}{\alpha^{2}}\left\{sT_{22}\left[(y_{k}b_{l}+y_{l}b_{k})\delta^{i}{}_{j}+a_{jl}b_{k}y^{i}\right]+\frac{1}{\alpha}(T-sT_{2}-s^{2}T_{22})
(y_{l}\delta^{i}{}_{j}+a_{lj}y^{i})y_{k}\right\}(k\rightarrow l\rightarrow j\rightarrow k)\nonumber\\
&&+\frac{c}{\alpha^{2}}\left[\frac{1}{\alpha^{3}}(3T-3sT_{2}-6s^{2}T_{22}-s^{3}T_{222})y_{k}y_{j}y_{l}+T_{222}b_{l}b_{k}b_{j}\right]y^{i}\nonumber\\
&&+\frac{c}{\alpha}\left[(H_{2}-sH_{22})(b_{j}-\frac{s}{\alpha}y_{j})a_{kl}-\frac{1}{\alpha^{2}}(H_{2}-sH_{22}-s^{2}H_{222})b_{l}y_{j}y_{k}
-\frac{sH_{222}}{\alpha}b_{k}b_{l}y_{j}\right]b^{i}(k\rightarrow l\rightarrow j\rightarrow k)\nonumber\\
&&+\frac{c}{\alpha}\left[\frac{s}{\alpha^{3}}(3H_{2}-3sH_{22}-s^{2}H_{222})y_{j}y_{k}y_{l}+H_{222}b_{l}b_{k}b_{j}\right]b^{i},\label{DC}
\end{eqnarray}
where $y_{i}:=a_{ij}y^{j}$ and $b^{i}:=a^{ij}b_{j}$, $c=c(x)\not=0$ is a scalar function on $M$.
\begin{eqnarray}
T:&=&-\frac{1}{n+1}[2sH+(b^{2}-s^{2})H_{2}],\label{T}\\
H:&=&\frac{\ptt-2(\po-s\pot)}{2\big[\p-s\pt+(b^2-s^2)\ptt\big]}.\label{H}
\end{eqnarray}
\end{prop}
\begin{proof}
By (\ref{bcij}), we have
\begin{eqnarray}\label{tu}
r_{00}=c\a^2,r_0=c\b,r=cb^2,r^i=cb^i,s^i{}_0=0,s_0=0,s^i=0.
\end{eqnarray}
Substituting (\ref{tu}) into (\ref{Gi}) yields
\begin{eqnarray}
G^i&=&{}^\a G^i+c\a\left\{\Theta(1+2Rb^2)+s\Omega\right\}y^i+c\a^2\left\{\Psi(1+2Rb^2)+s\Pi-R\right\}b^i\nonumber\\
&=&{}^\a G^i+c\a E y^i+c\a^2H b^i,\label{gp}
\end{eqnarray}
where
$$ E:=\frac{\pt+2s\po}{2\p}
-H\frac{s\p+(b^2-s^2)\pt}
{\p}.$$
Note that
\begin{eqnarray}\label{aby}
\alpha_{y^{i}}=\frac{y_{i}}{\alpha},~~ s_{y^{i}}=\frac{\alpha b_{i}-sy_{i}}{\a^{2}},
\end{eqnarray}
where $y_{i}:=a_{ij}y^{j}$.
\begin{eqnarray}
\frac{\partial G^{m}}{\partial y^{m}}=\frac{\partial {}^\a G^{m}}{\partial y^{m}}+c\a[(n+1)E+2sH+(b^{2}-s^{2})H_{2}],\label{Gmm}
\end{eqnarray}
where we take Einstein summation convention.
By (\ref{gp}) and (\ref{Gmm}), we have
\begin{eqnarray}\label{GFa}
G^{i}-\frac{1}{n+1}\frac{\partial G^{m}}{\partial y^{m}}y^{i}={}^\a G^{i}-\frac{1}{n+1}\frac{\partial {}^\a G^{m}}{\partial y^{m}}y^{i}+c\a(Ty^{i}+\alpha H b^{i}).
\end{eqnarray}
Put
\begin{eqnarray}\label{Wi}
W^{i}:=\a Ty^{i}+\alpha^{2} H b^{i}.
\end{eqnarray}
Differentiating (\ref{Wi}) with respect to $y^{j}$ yields
\begin{eqnarray}\label{Wij}
\frac{\partial W^{i}}{\partial y^{j}}=\a T\delta^{i}{}_{j}+(T\a_{y^{j}}+\alpha T_{2}s_{y^{j}})y^{i}+\big\{[\alpha^{2}]_{y^{j}} H+\a^{2}H_{2}s_{y^{j}}\big\} b^{i}.
\end{eqnarray}
Differentiating (\ref{Wij}) with respect to $y^{k}$ yields
\begin{eqnarray}\label{Wijk11}
\frac{\partial^{2} W^{i}}{\partial y^{j}\partial y^{k}}&=&\big[(T\a_{y^{k}}+\a T_{2}s_{y^{k}})\delta^{i}{}_{j}+T_{2}s_{y^{k}}\a_{y^{j}}y^{i}+H_{2}[\a^{2}]_{y^{j}}s_{y^{k}}b^{i}\big](k\leftrightarrow j)\nonumber\\
&&+\big(T\a_{y^{j}y^{k}}+\a T_{22}s_{y^{k}}s_{y^{j}}+\a T_{2}s_{y^{j}y^{k}}\big)y^{i}\nonumber\\
&&+\big\{[\a^{2}]_{y^{j}y^{k}}H+\a^{2}H_{22}s_{y^{k}}s_{y^{j}}+\a^{2}H_{2}s_{y^{j}y^{k}}\big\}b^{i},
\end{eqnarray}
where $k\leftrightarrow j$ denotes symmetrization. Therefore, it follows from (\ref{Wijk11}) that
\begin{eqnarray}\label{Wijk}
\frac{\partial^{3} W^{i}}{\partial y^{j}\partial y^{k}\partial y^{l}}&=& \big[T_{2}(\a_{y^{k}}s_{y^{l}}+\a_{y^{l}}s_{y^{k}}+\a s_{y^{k}y^{l}})+T\a_{y^{k}y^{l}}+\a T_{22}s_{y^{l}}s_{y^{k}}\big]\delta^{i}{}_{j}(k\rightarrow l\rightarrow j\rightarrow k)\nonumber\\
&&+\big[T_{2}(s_{y^{k}}\a_{y^{j}y^{l}}+\a_{y^{k}}s_{y^{j}y^{l}})+T_{22}(\a_{y^{k}}s_{y^{j}}+\a s_{y^{k}y^{j}})s_{y^{l}}\big]y^{i}
(k\rightarrow l\rightarrow j\rightarrow k)\nonumber\\
&&+\Large\{H_{2}\big([\a^{2}]_{y^{k}y^{l}}s_{y^{j}}+[\a^{2}]_{y^{k}}s_{y^{j}y^{l}}\big)+H_{22}\big([\a^{2}]_{y^{k}}s_{y^{l}}s_{y^{j}}
+\a^{2}s_{y^{k}y^{l}}s_{y^{j}}\big)\Large\}b^{i}(k\rightarrow l\rightarrow j\rightarrow k)\nonumber\\
&&+\big(T\a_{y^{j}y^{k}y^{l}}+\a T_{222}s_{y^{j}}s_{y^{k}}s_{y^{l}}+\a T_{2}s_{y^{j}y^{k}y^{l}}\big)y^{i}\nonumber\\
&&+\big\{H[\a^{2}]_{y^{j}y^{k}y^{l}}+\a^{2} H_{222}s_{y^{j}}s_{y^{k}}s_{y^{l}}+\a^{2} H_{2}s_{y^{j}y^{k}y^{l}}\big\}b^{i},
\end{eqnarray}
where $k\rightarrow l\rightarrow j\rightarrow k$ denotes cyclic permutation.
It follows from (\ref{aby}) that
\begin{eqnarray}\label{asd}
[\a^{2}]_{y^{l}}&=& 2y_{l},~~~ [\a^{2}]_{y^{l}y^{j}}=2a_{lj},~~~ [\a^{2}]_{y^{l}y^{j}y^{k}}=0,\label{a2d} \\
\a_{y^{l}y^{j}}&=& \frac{1}{\a}\big(a_{lj}-\frac{y_{l}}{\a}\frac{y_{j}}{\a}\big),~~ \a_{y^{l}y^{j}y^{k}}=-\frac{1}{\a^{3}}[a_{kl}y_{j}(k\rightarrow l\rightarrow j\rightarrow k)-\frac{3}{\a^{2}}y_{l}y_{j}y_{k}],\label{ad}\\
s_{y^{l}y^{j}}&=& -\frac{1}{\a^{2}}\big[s a_{lj}+\frac{1}{\a}(b_{l}y_{j}+b_{j}y_{l})-\frac{3s}{\a^{2}}y_{l}y_{j}\big],\label{s2d}\\
s_{y^{l}y^{j}y^{k}}&=& \frac{1}{\a^{5}}\big\{[\a(3s y_{j}-\a b_{j})a_{lk}+3b_{k}y_{l}y_{j}](k\rightarrow l\rightarrow j\rightarrow k)-\frac{15s}{\a}y_{k}y_{l}y_{j}\big\}.\label{s3d}
\end{eqnarray}
Plugging (\ref{a2d}), (\ref{ad}), (\ref{s2d}) and (\ref{s3d}) into (\ref{Wijk}) yields
\begin{eqnarray}
\frac{\partial ^{3}W^{i}}{\partial y^{j}\partial y^{k}\partial y^{l}}&=&\frac{1}{\a}\left\{\big[(T-sT_{2})a_{kl}+T_{22}b_{l}b_{k}\big]\delta^{i}{}_{j}+
\frac{1}{\a^{2}}\big[\frac{s}{\a}(3T_{22}+sT_{222})y_{l}-(T_{22}+sT_{222})b_{l}\big]y_{j}b_{k}y^{i}\right\}(k\rightarrow l\rightarrow j\rightarrow k)\nonumber\\
&&-\frac{1}{\a^{2}}\left\{sT_{22}\big[(y_{k}b_{l}+y_{l}b_{k})\delta^{i}{}_{j}+a_{jl}b_{k}y^{i}\big]
+\frac{1}{\a}(T-sT_{2}-s^{2}T_{22})(y_{l}\delta^{i}_{j}+a_{jl}y^{i})y_{k}\right\}(k\rightarrow l\rightarrow j\rightarrow k)\nonumber\\
&&+\frac{1}{\alpha^{2}}\left[\frac{1}{\alpha^{3}}(3T-3sT_{2}-6s^{2}T_{22}-s^{3}T_{222})y_{k}y_{j}y_{l}+T_{222}b_{l}b_{k}b_{j}\right]y^{i}\nonumber\\
&&+\frac{1}{\alpha}\left[(H_{2}-sH_{22})(b_{j}-\frac{s}{\alpha}y_{j})a_{kl}-\frac{1}{\alpha^{2}}(H_{2}-sH_{22}-s^{2}H_{222})b_{l}y_{j}y_{k}
-\frac{sH_{222}}{\alpha}b_{k}b_{l}y_{j}\right]b^{i}(k\rightarrow l\rightarrow j\rightarrow k)\nonumber\\
&&+\frac{1}{\alpha}\left[\frac{s}{\alpha^{3}}(3H_{2}-3sH_{22}-s^{2}H_{222})y_{j}y_{k}y_{l}+H_{222}b_{l}b_{k}b_{j}\right]b^{i},\label{Wijk1}
\end{eqnarray}
It follows from ${}^\a G^i(x,y)=\frac{1}{2}\Gamma^{i}_{jk}(x)y^{j}y^{k}$ that
\begin{eqnarray}\label{ajlk}
\frac{\partial ^{3}}{\partial y^{j}\partial y^{k}\partial y^{l}}\big[{}^\a G^{i}-\frac{1}{n+1}\frac{\partial {}^\a G^{m}}{\partial y^{m}}y^{i}\big]=0
\end{eqnarray}
By (\ref{DT}), (\ref{GFa}), (\ref{Wi}), (\ref{Wijk1}) and (\ref{ajlk}), we obtain (\ref{DC}).

\end{proof}
\section{Proof of Theorem \ref{main1}}
In this section, we mainly prove Theorem \ref{main1}. Firstly, we give the following Lemma.

\begin{lem}\label{HE}
Suppose that $\b$ satisfies (\ref{bcij}), then a general $\ab$-metric $F=\gab$ is a Douglas metric if and only if
$H_{2}-sH_{22}=0$, where $H$ is given by (\ref{H}).
\end{lem}
\begin{proof}
Suppose that a general $\ab$-metric $F$ is a Douglas metric, then the Douglas curvature of $F$ vanishes, i.e, $D^{i}{}_{jkl}=0$. From (\ref{bcij}) and (\ref{DC}), it follows that both rational and irrational parts of $D^{i}{}_{jkl}=0$ should vanish, i.e,
\begin{eqnarray}
&&\a^{4}\big\{[(T-sT_{2})a_{kl}+T_{22}b_{l}b_{k}]\delta^{i}{}_{j}+(H_{2}-sH_{22})a_{kl}b_{j}b^{i}+\frac{1}{3}H_{222}b_{l}b_{k}b_{j}b^{i}\big\}(k\rightarrow l\rightarrow j\rightarrow k)\nonumber\\
&&-\a^{2}\big[(T_{22}+sT_{222})b_{l}b_{k}y_{j}y^{i}+(T-sT_{2}-s^{2}T_{22})(y_{l}\delta^{i}{}_{j}+a_{jl}y^{i})y_{k}
+(H_{2}-sH_{22}-s^{2}H_{222})b_{l}y_{j}y_{k}b^{i}\big](k\rightarrow l\rightarrow j\rightarrow k)\nonumber\\
&&+(3T-3sT_{2}-6s^{2}T_{22}-s^{3}T_{222})y_{k}y_{j}y_{l}b^{i}=0,\label{rational}\\
&&\a^{2}\big\{\frac{1}{3}T_{222}b_{l}b_{k}b_{j}y^{i}-sT_{22}[(y_{k}b_{l}+y_{l}b_{k})\delta^{i}{}_{j}+a_{jl}b_{k}y^{i}]
-s[(H_{2}-sH_{22})a_{kl}y_{j}+H_{222}b_{k}b_{l}y_{j}]b^{i}\big\}(k\rightarrow l\rightarrow j\rightarrow k)\nonumber\\
&&+s(3T_{22}+sT_{222})y_{l}y_{j}b_{k}y^{i}(k\rightarrow l\rightarrow j\rightarrow k)+s(3H_{2}-3sH_{22}-s^{2}H_{222})y_{j}y_{k}y_{l}b^{i}=0,\label{irrational}
\end{eqnarray}
where $H$ and $T$ are given by (\ref{H}) and (\ref{T}), respectively. For $s\neq 0$, multiplying (\ref{irrational}) by $y^{j}y^{k}y^{l}$ yields
\begin{eqnarray}\label{TH2}
T_{222}y^{i}-\a H_{222}b^{i}=0.
\end{eqnarray}
By (\ref{TH2}),  it is easy to see that
\begin{eqnarray}\label{THe}
T_{222}=0, ~~ H_{222}=0.
\end{eqnarray}
Inserting (\ref{THe}) into (\ref{irrational}) yields
\begin{eqnarray}\label{irrational1}
&&\a^{2}s\big\{T_{22}[(y_{k}b_{l}+y_{l}b_{k})\delta^{i}{}_{j}+a_{jl}b_{k}y^{i}]+(H_{2}-sH_{22})a_{kl}y_{j}b^{i}\big\}(k\rightarrow l\rightarrow j\rightarrow k)\nonumber\\
&&-3sT_{22}y_{l}y_{j}b_{k}y^{i}(k\rightarrow l\rightarrow j\rightarrow k)-3s(H_{2}-sH_{22})y_{j}y_{k}y_{l}b^{i}=0.
\end{eqnarray}
Multiplying (\ref{irrational1}) by $b^{j}b^{k}b^{l}$ yields
\begin{eqnarray}\label{irrationa2}
b^{2}(b^{2}-3s^{2})T_{22}y^{i}+\a
s[2b^{2}T_{22}+(b^{2}-s^{2})(H_{2}-sH_{22})]b^{i}=0.
\end{eqnarray}
It follows from (\ref{irrationa2}) that
\begin{eqnarray}\label{THe1}
T_{22}=0, ~~~ H_{2}-sH_{22}=0.
\end{eqnarray}
Plugging (\ref{THe}) and (\ref{THe1}) into (\ref{rational}), we have
\begin{eqnarray}\label{rational1}
(T-sT_{2})\big\{\a^{2}[\a^{2}a_{kl}\delta^{i}{}_{j}-(y_{l}\delta^{i}{}_{j}+a_{jl}y^{i})y_{k}](k\rightarrow l\rightarrow j\rightarrow k)+3y_{k}y_{j}y_{l}y^{i}\big\}=0
\end{eqnarray}
Multiplying (\ref{rational1}) by $b^{j}b^{k}b^{l}$ yields
\begin{eqnarray}\label{rational2}
(T-sT_{2})(b^{2}-s^{2})(\a b^{i}-s y^{i})=0.
\end{eqnarray}
By (\ref{rational2}), we have
\begin{eqnarray}\label{TE}
T-sT_{2}=0.
\end{eqnarray}
By (\ref{T}), we obtain
\begin{eqnarray}\label{TEH}
T-sT_{2}=-\frac{1}{n+1}(b^{2}-s^{2})(H_{2}-sH_{22}).
\end{eqnarray}
By (\ref{TEH}), it is easy to see that the second equality of (\ref{THe1}) implies (\ref{TE}).

Conversely, suppose that the second equality of (\ref{THe1}) holds, it follows from (\ref{TEH}) that (\ref{TE}) holds.
Moreover,
\begin{eqnarray}\label{TH22}
T_{22}=0, ~~ T_{222}=0, ~~ H_{222}=0.
\end{eqnarray}
Plugging the second equality of (\ref{THe1}), (\ref{TE}) and (\ref{TH22}) into (\ref{DC}), we have $D^{i}{}_{jkl}=0$.
Hence, general $\ab$-metric $F=\gab$ is a Douglas metric.
\end{proof}

\begin{proof}[Proof of Theorem \ref{main1}]
By Lemma \ref{HE}, we obtain
\begin{eqnarray}\label{Hfg}
H=\frac{1}{2}[f(b^2)+g(b^2)s^2],
\end{eqnarray}
where $f$ and $g$ are two arbitrary smooth functions of $b^{2}$. By (\ref{H}) and (\ref{Hfg}), we will complete the proof of Theorem \ref{main1}.
\end{proof}
By taking $f=0$ and $g=0$, we obtain the following result \cite{szm-yct-oefm}
\begin{cor}
Let $F=\gab$ be a Finsler metric. Suppose that $\b$ satisfies (\ref{bcij}).
Then $F$ is projectively equivalent to $\alpha$ if and only if $\p(b^2,s)$ satisfies
\begin{eqnarray*}
\phi_{22}- 2 (\phi_1-s \phi_{12} ) =0.
\end{eqnarray*}
\end{cor}
\section{General $\ab$-metrics with vanishing Douglas curvature}
\begin{proof}[Proof of Theorem \ref{main2}]
Note that $(\phi-s\phi_{2})_{2}=-s\phi_{22}$. Therefore, (\ref{pde1****}) is changed to the following form
\begin{eqnarray}\label{PDE1}
2(\phi-s\phi_{2})_{1}+\frac{1}{s}[1-(f+gs^{2})(b^{2}-s^{2})](\phi-s\phi_{2})_{2}+(f+gs^{2})(\phi-s\phi_{2})=0.
\end{eqnarray}
Put
\begin{eqnarray}\label{psi}
\psi:=(\phi-s\phi_{2})\sqrt{b^{2}-s^{2}}.
\end{eqnarray}
Then
\begin{eqnarray}
\psi_{1}&=&(\phi-s\phi_{2})_{1}\sqrt{b^{2}-s^{2}}+\frac{1}{2\sqrt{b^{2}-s^{2}}}(\phi-s\phi_{2}), \label{phi11}\\ \psi_{2}&=&(\phi-s\phi_{2})_{2}\sqrt{b^{2}-s^{2}}-\frac{s}{\sqrt{b^{2}-s^{2}}}(\phi-s\phi_{2}). \label{phi2}
\end{eqnarray}
It follows from (\ref{phi11}) and (\ref{phi2}) that (\ref{PDE1}) is equivalent to
\begin{eqnarray}\label{PDE2}
\psi_{1}+\frac{1}{2s}\big[1-(f+gs^{2})(b^{2}-s^{2})\big]\psi_{2}=0.
\end{eqnarray}
The characteristic equation of PDE (\ref{PDE2}) is
\begin{eqnarray}\label{characteristicE}
\frac{db^{2}}{1}=\frac{ds}{\frac{1}{2s}\big[1-(f+gs^{2})(b^{2}-s^{2})\big]}
\end{eqnarray}
(\ref{characteristicE}) is equivalent to
\begin{eqnarray}\label{equivE}
2s\frac{ds}{db^{2}}=1-(f+gs^{2})(b^{2}-s^{2}).
\end{eqnarray}
Set
\begin{eqnarray}\label{chi}
\chi(b^{2})=s^{2}(b^{2})-b^{2}.
\end{eqnarray}
Plugging (\ref{chi}) into (\ref{equivE}) yields
\begin{eqnarray*}\label{Bernoulli}
\frac{d\chi}{db^{2}}=(f+gb^{2})\chi+g\chi^{2}.
\end{eqnarray*}
This is a Bernoulli equation which can be rewritten as
 \begin{eqnarray*}
 \frac{d}{db^{2}}\big(\frac{1}{\chi}\big)=-(f+gb^{2})\frac{1}{\chi}-g.
 \end{eqnarray*}
This is a linear 1-order ODE of $\frac{1}{\chi}$.
One can easily get its solution
 \begin{eqnarray}\label{solution}
 \frac{1}{\chi}=-e^{-\int(f+gb^{2})db^{2}}\big[c+\int g e^{\int(f+gb^{2})db^{2}}db^{2}],
  \end{eqnarray}
 where $c$ is an arbitrary constant.
By (\ref{chi}) and (\ref{solution}), the independent integral of (\ref{characteristicE}) is
\begin{eqnarray*}
\frac{b^{2}-s^{2}}{e^{\int(f+gb^{2})db^{2}}-(b^{2}-s^{2})\int g e^{\int(f+gb^{2})db^{2}}db^{2}}=\frac{1}{c}.
\end{eqnarray*}
Hence the solution of (\ref{PDE2}) is
\begin{eqnarray}\label{solution1}
\psi=\Phi\left(\frac{b^{2}-s^{2}}{e^{\int(f+gb^{2})db^{2}}-(b^{2}-s^{2})\int g e^{\int(f+gb^{2})db^{2}}db^{2}}\right),
\end{eqnarray}
where $\Phi$ is any continuously differentiable function.
By (\ref{psi}) and (\ref{solution1}), we have
\begin{eqnarray}\label{phi1}
\phi-s\phi_{2}=\Phi\left(\frac{b^{2}-s^{2}}{e^{\int(f+gb^{2})db^{2}}-(b^{2}-s^{2})\int g e^{\int(f+gb^{2})db^{2}}db^{2}}\right)\frac{1}{\sqrt{b^{2}-s^{2}}}.
\end{eqnarray}
Let $\phi=s\varphi$, then we have
\begin{eqnarray}\label{phis2}
\phi-s\phi_{2}=-s^{2}\varphi_{2}.
\end{eqnarray}
By (\ref{phi1}) and (\ref{phis2}), we obtain
\begin{eqnarray}
\varphi=h(b^{2})-\int \frac{\Phi\big(\eta(b^{2},s)\big)}{s^{2}\sqrt{b^{2}-s^{2}}}ds,
\end{eqnarray}
where $h(x)$ is an arbitrary smooth function and $\eta(b^{2},s)$ is given by (\ref{eta}).  Hence, by $\phi=s\varphi$, we get (\ref{mainsoltion}).
\end{proof}
In the following, we will give necessary and sufficient conditions for a general $\ab$-metric with vanishing Douglas curvature to be a Finsler metric.
\begin{lem}\label{finsler}
 Let $F=\gab$ be a general $\ab$-metric on an $n$-dimensional manifold $M$, where $\phi$ is given by (\ref{mainsoltion}). Then $F$ is a Finsler metric if and only if
  \begin{eqnarray}\label{ineq1}
  \frac{\Phi}{\sqrt{b^{2}-s^{2}}}>0, ~~ -\frac{\sqrt{b^{2}-s^{2}}}{s}\Phi_{2}>0.
  \end{eqnarray}
  when $n\geq 3$ or
  \begin{eqnarray}\label{ineq2}
   -\frac{\sqrt{b^{2}-s^{2}}}{s}\Phi_{2}>0.
  \end{eqnarray}
  when $n=2$.
\end{lem}
\begin{proof}
Note that $-s\phi_{22}=(\phi-s\phi)_{2}$. By (\ref{ppp}), (\ref{ppp1}) and (\ref{phi1}), we will get (\ref{ineq1}) and (\ref{ineq2}).
\end{proof}
\section{Some new examples}
In this section, we will explicitly construct some new examples .
\begin{example}
Take $g=0$ and $\Phi(\eta(b^{2},s))=(b^{2}-s^{2})^{\frac{m}{2}}e^{-\int fdb^{2}}$, then for any natural number $m\geq 1$, parts of the solutions of
(\ref{mainsoltion}) are given by
\begin{enumerate}[(a)]
\item $m=2l$
\begin{eqnarray*}
 \phi(b^{2},s)&=& \tilde{h}_{1}(b^{2})s-e^{-\int f db^{2}}\frac{(2l-1)!!}{(2l-2)!!}\Large\{\sum_{i=1}^{l-1}\frac{(2l-2-2i)!!}{(2l-2i+1)!!}
 \big(b^{2}\big)^{i-1}(b^{2}-s^{2})^{\frac{2l-2i+1}{2}}\nonumber\\
 &&- \big(b^{2}\big)^{l-1}\big[(b^{2}-s^{2})^{\frac{1}{2}}+s \arctan\frac{s}{\sqrt{b^{2}-s^{2}}}\big]\Large\},
 \end{eqnarray*}
\item $m=2l+1$
\begin{eqnarray*}
\phi(b^{2},s)&=& \tilde{h}_{2}(b^{2})s-e^{-\int f db^{2}} \frac{(2l)!!}{(2l-1)!!}\left[\sum_{i=1}^{l}\frac{(2l-2i-1)!!}{(2l-2i+2)!!}\big(b^{2}\big)^{i-1}(b^{2}-s^{2})^{l-i+1}-\big(b^{2}\big)^{l}\right],
\end{eqnarray*}
\end{enumerate}
where $\tilde{h}_{1}$, $\tilde{h}_{2}$ and $f$ are any smooth functions of $b^{2}$ such that $\phi$ is positive. Moreover, the corresponding general $\ab$-metrics
$$
F=\gab
$$
are of Douglas type.
\end{example}

Note that we have made use of Lemma \ref{formu}.
\begin{example}
Take $g=0$, $f=\frac{\mu^{2}+\varepsilon\xi}{\varepsilon+(\mu^{2}+\varepsilon\xi)b^{2}}$ and $\Phi(\eta(b^{2},s))=\varepsilon\sqrt{\frac{\eta}{1-\mu^{2}\eta}}$, then parts of the solutions of
(\ref{mainsoltion}) are given by
\begin{eqnarray}\label{example2}
\phi(b^{2},s)=\tilde{h}(b^{2})s+\frac{\sqrt{\varepsilon+\varepsilon\xi b^{2}+\mu^{2}s^{2}}}{1+\xi b^{2}},
\end{eqnarray}
where $\tilde{h}$ is a smooth function of $b^{2}$ and $\mu$, $\varepsilon$, $\xi$ are constants such that $\phi$ is positive. Moreover, the corresponding general $\ab$-metrics
$$
F=\gab
$$
are of Douglas type.
\end{example}
\textbf{Remark:}
Especially, take $\tilde{h}(b^{2})=\frac{\mu}{1+\xi b^{2}}$ in (\ref{example2}), we have
\begin{eqnarray}\label{funk}
\phi(b^{2},s)=\frac{\sqrt{\varepsilon+\varepsilon\xi b^{2}+\mu^{2}s^{2}}}{1+\xi b^{2}}+\frac{\mu s}{1+\xi b^{2}}.
\end{eqnarray}
\begin{enumerate}[(1)]
\item Take $\a=|y|$ and $\b=\langle x,y\rangle$, then the corresponding general $\ab$-metrics of (\ref{funk})
\begin{eqnarray*}
F=\frac{\sqrt{\varepsilon(1+\xi|x|^{2})+\mu^{2}\langle x,y\rangle^{2}}}{1+\xi|x|^{2}}+\frac{\mu\langle x,y\rangle^{2}}{1+\xi|x|}
\end{eqnarray*}
are of Douglas type. In fact, they belong to spherically symmetric Douglas metrics, too. Moreover, when $\varepsilon=1$, $\xi=-1$ and $\mu=\pm 1$, $F$ is just the Funk metric.
\item Take $\a=|y|$ and $\b=\langle x,y\rangle+\langle a,y\rangle$, where $a$ is a constant vector, then the corresponding general $\ab$-metrics of (\ref{funk})
    $$F=\f{\sqrt{(1-\xx-2\langle a,x\rangle-|a|^2)\yy+(\xy+\langle a,x\rangle)^2}}{1-\xx-2\langle a,x\rangle-|a|^2}
\pm\f{\xy+\langle a,x\rangle}{1-\xx-2\langle a,x\rangle-|a|^2}$$
are of Douglas type(See Example 8.1 in \cite{yu-zhu}). Actually, they are just the generalized Funk metrics expressed in some other local coordinate system.
\end{enumerate}
\begin{example}
Take $f=g=0$ and $\Phi(\eta(b^{2},s))=(1+\eta)\sqrt{\eta}$, then parts of the solutions of
(\ref{mainsoltion}) are given by
\begin{eqnarray}\label{example3}
\phi(b^{2},s)=\tilde{h}(b^{2})s+1+b^{2}+s^{2},
\end{eqnarray}
where $\tilde{h}$ is a smooth function of $b^{2}$ such that $\phi$ is positive. Moreover, the corresponding general $\ab$-metrics
$$
F=\gab
$$
are of Douglas type.
\end{example}

\textbf{Remark:} Take $\tilde{h}(b^{2})=2\sqrt{1+b^{2}}$ in (\ref{example3}), $\a=\f{\sqrt{(1+\mu|x|^2)|y|^2-\mu\langle x,y\rangle^2}}{1+\mu|x|^2}$ and $\b=\frac{\langle x,y\rangle}{(1+\mu|x|^2)\frac{3}{2}}$, where $\mu$ is a constant. We obtain Example 4.3 given in \cite{yct-zhm-onan}, namely
$$
F=\frac{(\sqrt{1+(1+\mu)|x|^{2}}\sqrt{(1+\mu|x|^{2})|y|^{2}-\mu\langle x,y\rangle^{2}}+\langle x,y\rangle)^{2}}{(1+\mu|x|^{2})^{2}\sqrt{(1+\mu|x|^{2})|y|^{2}-\mu\langle x,y\rangle^{2}}}.
$$
In particular, $F$ is just the Berwald's metric when $\mu=-1$.
\begin{example}
Take $f=g=0$ and $\Phi(\eta(b^{2},s))=\frac{\sqrt{\eta}}{(1-\eta)^{\frac{3}{2}}}$, then parts of the solutions of
(\ref{mainsoltion}) are given by
\begin{eqnarray}\label{example4}
\phi(b^{2},s)=\tilde{h}(b^{2})s+\frac{1-b^{2}+2s^{2}}{(1-b^{2})\sqrt{1-b^{2}+s^{2}}},
\end{eqnarray}
where $\tilde{h}$ is a smooth function of $b^{2}$ such that $\phi$ is positive. Moreover, the corresponding general $\ab$-metrics
$$
F=\gab
$$
are of Douglas type.
\end{example}
\textbf{Remark:} Take $\tilde{h}(b^{2})=\mp\frac{2}{(1-b^{2})^{2}}$ in (\ref{example4}), $\a=|y|$ and $\b=\langle x,y\rangle+\langle a,y\rangle$, where $a$ is a constant vector, the corresponding general $\ab$-metrics of (\ref{example4})
$$F=\frac{\{\sqrt{(1-\xx-2\langle a,x\rangle-|a|^2)\yy+(\xy+\langle a,y\rangle)^2}\mp(\xy+\langle a,y\rangle)\}^2}{(1-\xx-2\langle a,x\rangle-|a|^2)^2\sqrt{(1-\xx-2\langle a,x\rangle-|a|^2)\yy+(\xy+\langle a,y\rangle)^2}}$$
are of Douglas type(See Example 8.2 in \cite{yu-zhu}). Actually, they are just the generalized Berwald's metrics expressed in some other local coordinate system.
\begin{example}
Take $f=g=0$ and $\Phi(\eta(b^{2},s))=\frac{1}{2}[\frac{1}{\sqrt{c-\eta}}-\frac{\varepsilon}{\sqrt{c-\varepsilon^{2}\eta}}]\sqrt{\eta}$, then parts of the solutions of
(\ref{mainsoltion}) are given by
\begin{eqnarray}\label{example5}
\phi(b^{2},s)=\tilde{h}(b^{2})s+\frac{1}{2}\left[\frac{\sqrt{c-b^{2}+s^{2}}}{c-b^{2}}
-\frac{\varepsilon\sqrt{c-\varepsilon^{2}(b^{2}-s^{2})}}{c-\varepsilon^{2}b^{2}}\right],
\end{eqnarray}
where $c>0$, $\varepsilon<1$ , $\tilde{h}$ is a smooth function of $b^{2}$ such that $\phi$ is positive. Moreover, the corresponding general $\ab$-metrics
$$
F=\gab
$$
are of Douglas type.
\end{example}
\textbf{Remark:} Take $\tilde{h}(b^{2})=\frac{1}{2}\left(\frac{1}{c-b^{2}}-\frac{\varepsilon^{2}}{c-\varepsilon^{2}b^{2}}\right)$ in (\ref{example5}),
$\a=|y|$ and $\b=\langle x,y\rangle+\langle a,y\rangle$, where $a$ is a constant vector, the corresponding general $\ab$-metrics of (\ref{example5})
\begin{eqnarray*}
F&=&\f{1}{2}\Bigg\{\frac{\sqrt{(c-\xx-2\langle a,x\rangle-|a|^2)\yy+(\xy+\langle a,y\rangle)^2}+\xy+\langle a,y\rangle}{c-\xx-2\langle a,x\rangle-|a|^2}\\
&&-\frac{\varepsilon\sqrt{[c-\varepsilon^2(\xx+2\langle a,x\rangle+|a|^2)]\yy+\varepsilon^2(\xy+\langle a,y\rangle)^2}+\varepsilon^2(\xy+\langle a,y\rangle)}{c-\varepsilon^2(\xx+2\langle a,x\rangle+|a|^2)}\Bigg\}
\end{eqnarray*}
are of Douglas type. In particular, when $c=1$ and $a=0$, it is just Shen' metrics(see (39) in \cite{szm-proj}). When $c=1$, it is just the Example 8.4 in \cite{yu-zhu}. When $a=0$, it is just a projectively flat spherically symmetric Finsler metrics with constant flag curvature $-1$ \cite{Mo-Zhu}.
\begin{example}
Take $f=\lambda$, $g=\frac{\lambda^{2}}{1-\lambda b^{2}}$ and $\Phi(\eta(b^{2},s))=\sqrt{\eta}$, then parts of the solutions of
(\ref{mainsoltion}) are given by
\begin{eqnarray}\label{example6}
\phi(b^{2},s)=\tilde{h}(b^{2})s+\frac{\sqrt{(1-\lambda b^{2})(1-2\lambda b^{2}+\lambda s^{2})}}{1-2\lambda b^{2}},
\end{eqnarray}
where $\lambda$ is an arbitrary constant, $\tilde{h}$ is a smooth function of $b^{2}$ such that $\phi$ is positive. Moreover, the corresponding general $\ab$-metrics
$$
F=\gab
$$
are of Douglas type.
\end{example}
\textbf{Remark:} Take $\tilde{h}(b^{2})=\frac{\sqrt{1-\lambda b^{2}}}{1-2\lambda b^{2}}$ in (\ref{example6}),
$\a=|y|$ and $\b=\langle x,y\rangle+\langle a,y\rangle$, where $a$ is a constant vector, the corresponding general $\ab$-metrics of (\ref{example6})
\begin{eqnarray*}
F=\frac{\sqrt{1-\lambda(\xx+2\langle a,x\rangle+|a|^2)}}{1-2\lambda(\xx+2\langle a,x\rangle+|a|^2)}\{\sqrt{(1-2\lambda(\xx+2\langle a,x\rangle+|a|^2))|y|^{2}+\lambda(\langle x,y\rangle+\langle a,y\rangle)^{2}}+\langle x,y\rangle+\langle a,y\rangle\}
\end{eqnarray*}
are of Douglas type, but not locally projectively flat.

\noindent Hongmei Zhu\\
College of Mathematics and Information Science, Henan Normal University, Xinxiang, 453007, P.R. China\\
zhm403@163.com
\end{document}